\documentclass{elsarticle}
\usepackage{amsmath,amssymb,latexsym}
\input{xypic}

\newtheorem{thm}{Theorem}
\newtheorem{lem}[thm]{Lemma}
\newtheorem{prop}[thm]{Proposition}
\newtheorem{cor}[thm]{Corollary}
\newdefinition{dfn}{Definition}

\newdefinition{rmk}{Remark}
\newdefinition{ex}{Example}
\newdefinition{ntn}{Notations}

\newproof{pf}{proof}

\numberwithin{equation}{section}

\begin{document}

\hyphenation{con-si-de-ring o-pe-ra-tion mul-ti-o-pe-ra-tions  mul-ti-o-pe-ra-tion
ar-bi-tra-ry ge-ne-ra-li-za-tion ge-ne-ra-li-za-tions cha-rac-te-ri-za-tion to-pics
sub-ca-te-go-ry ca-te-go-ry cor-res-pon-ding hy-per-groups}

\begin{frontmatter}

\title{Fundamental relations in multialgebras. Applications}

\author{Cosmin Pelea\fnref{1,2}}

\address{Babe\c s-Bolyai University, Faculty of Mathematics and
Computer Science, M. Kog\u alniceanu 1, RO-400084
Cluj-Napoca, Romania} \ead{cpelea@math.ubbcluj.ro}
\fntext[1]{Corresponding author.}
\fntext[2]{The first author was supported by the CNCS-UEFISCDI grant
PN-II-RU-TE-2011-3-0065.}

\author{Ioan Purdea}

\address{Babe\c s-Bolyai University, Faculty of Mathematics and
Computer Science, M. Kog\u alniceanu 1, RO-400084
Cluj-Napoca, Romania} \ead{purdea@math.ubbcluj.ro}

\author{Liana Stanca}

\address{Babe\c s-Bolyai University, Faculty of Economic Sciences and Business Administration,
Business Information Systems Dpt., T. Mihali 58-60, RO-400591
Cluj-Napoca, Romania} \ead{liana.stanca@econ.ubbcluj.ro}


\begin{abstract}
Based on the properties of the poset of those equivalence relations of a multialgebra
for which the factor multialgebra is a universal algebra, we give a characterization for the
fundamental relations of a multialgebra. We point out the benefits of our approach
by giving two applications. One of them provides a new characterization of the commutative
fundamental relation of a hyperring, and the other will give a general category theoretical
property of the construction of the fundamental algebras (both in the general case and in the 
hyperring case).
\end{abstract}

 \end{frontmatter}

\section{Introduction}

Multialgebras (also called hyperstructures) appeared more than 70 years ago and in the last 30 years one can notice 
an increased interest in them as well as many attempts to connect them to other mathematical
or computer science topics (see, for instance, \cite{CoLe}). In the 60's the papers \cite{Gr62} and 
\cite{Pi67} were published and showed that maybe the most natural connection of multialgebras 
is the one with universal algebras. Many computer science theorists use multialgebras
to model (arbitrary) relational systems. Our multialgebras are the structures characterized in \cite{Gr62}. From our point 
of view, arbitrary relational systems can be modeled only as partial multialgebras. There are many differences  
between the two clases of algebraic structures. Some of them can be found in \cite{Pe04}.

The first part of our paper refers to the connection between multialgebras and universal algebras 
and investigates those equivalence relations of a multialgebra for which the factor multialgebra is a universal algebra.
These relations are generalizations of universal algebras verbal congruences (see \cite{BeBe}) and,
lately, they can be found in the literature under the generic name ``{\it fundamental relations}'', while 
the corresponding factor (multi)algebras are named ``{\it fundamental algebras}`` (see, for instance, \cite{MiDa}).
The first important steps in characterizing these relations for general multialgebras can be found 
in \cite{BrPe} and \cite{Pe01}. 
Using some results from \cite{PePu06} and some properties of the complete lattice of the
equivalence relations of a multialgebra $\mathfrak A$ for which the factor multialgebra is a 
universal algebra, we give a characterization for the smallest equivalence relation of a
multialgebra for which the factor multialgebra is a universal algebra satisfying a set $\mathcal I$ of identities.
Using the example given by \cite{MiDa}, we propose to name these relations $\mathcal I$-{\it fundamental relations} for a
better distinction between them and the relation $\alpha^*_{\mathfrak A}$ (see Section \ref{s3})
which is known as {\it the fundamental relation of} $\mathfrak A$.
The importance of the results from Section \ref{s3} for hyperstructure theory consists 
mainly in the fact that they give a general ``recipe'',
but also in the way they are proved without developing a large number of computations. 
It must be said that a previous version of the proof of Theorem \ref{qrI} 
involved some unary polynomial functions computations. We are grateful to an anonymous referee of a 
preprint of our paper for the suggestions which lead to the current version of its proof. 

We connect our results with \cite{DaVo}, but they can also be used for many other 
hyperstructures. The advantage of our approach can be seen in the applications given in
Section \ref{sectionhr} and Section \ref{dirlim}. In Section \ref{sectionhr} we show
how our general results can be used for hyperrings for obtaining some of the results from \cite{DaVo}
using less computations. We also give in Theorem \ref{t3} a new characterization of the 
smallest equivalence relation of a hyperring for which the factor is a commutative ring 
(called $\alpha^*$-{\it relation} in \cite{DaVo} or {\it commutative fundamental relation} in \cite{MiDa}).   
In Section \ref{dirlim} we show that the construction of the $\mathcal I$-fundamental algebra
of a multialgebra determines a functor which preserves the directed colimits of multialgebras, and that 
this also works for hyperrings and their commutative fundamental rings. For the category
theoretical background of this section, see \cite{AHS}.

\section{Preliminaries}

Let  $\tau=(n_{\gamma})_{\gamma<o(\tau)}$ be a sequence of
nonnegative integers ($o(\tau)$ is an ordinal) and for any
$\gamma<o(\tau)$, let ${\mathbf f}_\gamma$ be a symbol of an
$n_{\gamma}$-ary (multi)operation. 
If $A$ is a set, we denote by $P^*(A)$ the set of the nonempty subsets of $A$.
A {\it multialgebra} $\mathfrak A$ {\it of type $\tau$}
consists in a set $A$ and a family of multioperations
$(f_{\gamma})_{\gamma < o(\tau)}$, where $f_{\gamma}:
A^{n_{\gamma}}\rightarrow P^*(A)$ is the $n_{\gamma}$-ary
multioperation which corresponds to the symbol ${\mathbf
f}_\gamma$. It is easy to observe that universal algebras are particular multialgebras.

Each multialgebra $\mathfrak A$ determines a
universal algebra $\mathfrak P^*(\mathfrak A)$ on $P^*(A)$ defining (as in \cite{Pi67}) for any
$\gamma<o(\tau)$, and any $A_0,\ldots ,A_{n_{\gamma}-1} \in
P^*(A)$,
$$f_{\gamma}(A_0,\ldots ,A_{n_{\gamma}-1})=\bigcup
\{f_{\gamma}(a_0,\ldots ,a_{n_{\gamma}-1}) \mid \  a_i\in A_i,\ i
\in \{0,\ldots ,n_{\gamma}-1 \} \}.$$  We call $\mathfrak
P^*(\mathfrak A)$ {\it the (universal) algebra of the nonempty
subsets of} $\mathfrak A$.

Since $\mathfrak P^*(\mathfrak A)$ is a universal algebra, one can
construct the universal algebra $\left({\rm Pol}_n(\mathfrak P^*(\mathfrak A)),(f_{\gamma})_{\gamma < o(\tau)}\right)$
of the $n$-ary polynomial functions of $\mathfrak P^*(\mathfrak A)$ (see, for instance, \cite{BuSa}). 
We denote by ${\rm Pol}^A_n (\mathfrak P^*(\mathfrak A))$ its subuniverse generated by all the $n$-ary 
singleton-providing constant functions and all the $n$-ary projections, i.e. by 
$\{c^n_a\mid a\in A\}\cup\{e_i^n\mid i\in\{0,\ldots,n-1\}\},$
where $c^n_a,e^n_i:P^*(A)^n\rightarrow P^*(A)$, 
$$c^n_a(A_0,\ldots,A_{n-1})=\{a\}\ \hbox{and}\
e^n_i(A_0,\ldots,A_{n-1})=A_i.$$ Clearly, the set ${\rm Clo}_n(\mathfrak P^*(\mathfrak A))$ of the $n$-ary term
functions on $\mathfrak P^*(\mathfrak A)$ is the subuniverse of
$\left({\rm Pol}^A_n (\mathfrak P^*(\mathfrak A)),(f_{\gamma})_{\gamma < o(\tau)})\right)$ generated by
$\{e_i^n\mid i\in\{0,\ldots,n-1\}\}$, and each $n$-ary term $\mathbf p$ of type $\tau$ determines an 
$n$-ary term function (denoted by) $p$ from ${\rm Clo}_n(\mathfrak P^*(\mathfrak A))$.

A map $h:A\rightarrow B$ between the  multialgebras $\mathfrak A$
and $\mathfrak B$ of the same type $\tau$ is called {\it homomorphism}
if for any $\gamma < o(\tau)$ and $a_0,\ldots,a_{n_{\gamma}-1}\in
A$ we have
$$
h(f_{\gamma}(a_0,\ldots,a_{n_{\gamma}-1}))\subseteq
f_{\gamma}(h(a_0),\ldots,h(a_{n_{\gamma}-1})).\eqno(1)
$$
A {\it multialgebra isomorphism} is a bijective map $h$ such that 
for any $\gamma < o(\tau)$ and $a_0,\ldots,a_{n_{\gamma}-1}\in A$ we have
$$
h(f_{\gamma}(a_0,\ldots,a_{n_{\gamma}-1}))=
f_{\gamma}(h(a_0),\ldots,h(a_{n_{\gamma}-1})).\eqno(1')
$$

Let ${\mathbf q}, {\mathbf r}$ be two $n$-ary terms of type $\tau$. 
The $n$-ary ({\it strong}) {\it identity}
${\mathbf q} = {\mathbf r}$ is satisfied on a
multialgebra $\mathfrak A$ if
$$q(a_0,\ldots ,a_{n-1}) = r(a_0,\ldots ,a_{n-1}),\ \forall
a_0,\ldots ,a_{n-1}\in A.$$ We say that the {\it weak identity}
${\mathbf q} \cap {\mathbf r}\neq\emptyset $ is satisfied on a
multialgebra $\mathfrak A$ if
$$q(a_0,\ldots ,a_{n-1}) \cap r(a_0,\ldots ,a_{n-1}) \neq
\emptyset,\ \forall a_0,\ldots ,a_{n-1}\in A$$ ($q$ and $r$ denote
the term functions induced by ${\mathbf q}$ and $ {\mathbf r}$
respectively on $\mathfrak P^*(\mathfrak A)$).

Let $\mathfrak A$ be a multialgebra of type $\tau$ and let $\rho$
be an equivalence relation on $A$. We denote by
$\rho \langle x \rangle $ the class of $x $ modulo $\rho$, and
$A/\rho=\{\rho \langle x \rangle \mid x\in A\}.$ Taking 
$$f_{\gamma}(\rho \langle a_0
\rangle ,\ldots ,\rho \langle a_{n_{\gamma}-1} \rangle )=\{\rho
\langle b \rangle \mid b\in f_{\gamma}(b_0,\ldots
,b_{n_{\gamma}-1}),\, a_i\rho b_i,\, i=0,\ldots ,n_{\gamma}-1\},$$ for each
$\gamma<o(\tau)$, one obtains a multialgebra $\mathfrak A/\rho$ on $A/\rho$ called
{\it the factor multialgebra of $\mathfrak A$ determined by
$\rho$} (see \cite{Gr62}).

\begin{rmk}
If a weak (or strong) identity
${\mathbf q} \cap {\mathbf r}\neq\emptyset$ (or ${\mathbf q} = {\mathbf r}$)  
is satisfied in a multialgebra $\mathfrak A$,
the class $\rho \langle a \rangle$ of any common element 
$a$ of $q(a_0,\ldots ,a_{n-1})$ and $r(a_0,\ldots ,a_{n-1})$ is an element of  
$q(\rho \langle a_0
\rangle ,\ldots ,\rho \langle a_{n-1} \rangle)\cap 
r(\rho \langle a_0
\rangle ,\ldots ,\rho \langle a_{n-1} \rangle),$
hence the weak identity ${\mathbf q} \cap {\mathbf r}\neq\emptyset$ is satisfied in the factor multialgebra 
$\mathfrak A/\rho$.
\end{rmk}

Let $\rho$ be an equivalence relation on the set $A$. We denote by
$\overline{\overline{\rho}}$ the relation defined on $P^*(A)$ as
follows: if $X,Y\in P^*(A)$, then
$$
X\overline{\overline{\rho}}Y\ \Leftrightarrow\ x\rho y,\ \forall
x\in X,\ \forall y\in Y\ (\Leftrightarrow\ X\times
Y\subseteq\rho).
$$
Denote by $E_{ua}(\mathfrak A)$, or $E_{ua}(A,(f_\gamma)_{\gamma<o(\tau)})$,
the set of all equivalence relations $\rho$ of a multialgebra
$\mathfrak A=(A,(f_\gamma)_{\gamma<o(\tau)})$ for which $\mathfrak A/\rho$ is a
universal algebra.

\begin{prop}{\rm\cite[Proposition 4.1]{PePu06}}\label{eua}
Let $\mathfrak
A=(A,(f_\gamma)_{\gamma<o(\tau)})$ be a multialgebra and let
$\rho$ be an equivalence relation on $A$. The following conditions
are equivalent:
\begin{itemize}
\item[{\rm a)}] $\rho\in E_{ua}(\mathfrak A)$;
\item[{\rm b)}] If $\gamma<o(\tau)$, $a,b,x_0,\ldots,x_{n_\gamma-1}\in A$ and $a\rho b$ 
then, for all $i\in\{0,\ldots,n_\gamma-1\}$,
$$f_\gamma(x_0,\ldots ,x_{i-1},a,x_{i+1},\ldots ,x_{n_\gamma-1})
\overline{\overline{\rho}} f_\gamma(x_0,\ldots
,x_{i-1},b,x_{i+1},\ldots ,x_{n_\gamma-1});$$
\item[{\rm c)}] If $\gamma<o(\tau)$, $x_i,y_i\in A$ and $x_i\rho y_i$
for any $i\in \{0,\ldots ,n_\gamma-1\}$, then
$$
f_\gamma(x_0,\ldots,x_{n_\gamma-1}) \overline{\overline{\rho}}
f_\gamma(y_0,\ldots,y_{n_\gamma-1});
$$
\item[{\rm d)}] If $n\in \mathbb N,$ $p\in {\rm Pol}^A_n(\mathfrak
P^*(\mathfrak A))$, $x_i,y_i\in A$ and $x_i\rho y_i$ $(i\in
\{0,\ldots,n-1\})$, then
$$
p(x_0,\ldots,x_{n-1})\overline{\overline{\rho}}
p(y_0,\ldots,y_{n-1}).
$$
\end{itemize}
\end{prop}

It is easy to notice that 
$E_{ua}(A,(f_\gamma)_{\gamma<o(\tau)})=\bigcap_{\gamma<o(\tau)}E_{ua}(A,f_\gamma).$

\begin{ex}\label{congA}
If $\mathfrak A$ is a universal algebra, then $E_{ua}(\mathfrak A)$ is the set $C(\mathfrak A)$ of the congruence
relations of $\mathfrak A$.
\end{ex}

\begin{ex}
We remind that a multialgebra $(A,\circ)$ with one binary multioperation is called {\it hypergroupoid}.
An equivalence relation $\rho$ of $A$ is in $E_{ua}(A,\circ)$ if and only if for any $a,b,x\in A,$ with
$a\rho b$, we have
$$a\circ x\,\overline{\overline{\rho}}\, b\circ x,\ x\circ a\,\overline{\overline{\rho}}\, x\circ b.$$
In the literature, such a relation is called {\it strongly regular equivalence relation}. The same term is sometimes 
used for general multialgebras $\mathfrak A$ to name the relations from $E_{ua}(\mathfrak A)$. 
The relations from $E_{ua}(\mathfrak A)$ can also be found in the literature as {\it strong 
congruence relations} (see \cite{AmRo09}). 
\end{ex}

\begin{rmk}\label{r2}
If $\rho\in E_{ua}(\mathfrak A)$ then in $\mathfrak A/\rho$, 
which is a universal algebra, 
$$f_{\gamma}(\rho \langle a_0 \rangle ,\ldots ,\rho \langle
a_{n_{\gamma}-1} \rangle )=\rho\langle b\rangle,\ \forall\,
b\in f_{\gamma}(a_0,\ldots
,a_{n_{\gamma}-1}).$$ 
It follows that if $n\in\mathbb N$ and $\mathbf p$ is an $n$-ary term of type $\tau$, then
$$p(\rho \langle a_0 \rangle ,\ldots ,\rho \langle a_{n-1} \rangle
)=\rho\langle b\rangle,\ \forall\,
b\in p(a_0,\ldots ,a_{n-1})\eqno(2)$$
(see \cite[Remark 13]{PePu08}), hence any
(weak or strong) identity of $\mathfrak A$ is also
satisfied on the algebra $\mathfrak A/\rho$.
\end{rmk}

\begin{ex}\label{exhypslash}
An $H_v$-{\it semigroup} is a hypergroupoid with a weak associative multioperation. 
If the multioperation is associative, the multialgebra is called {\it semihypergroup}. 
An $H_v$-{\it group} is an $H_v$-semigroup $(A,\circ)$ which satisfies the condition
$$a\circ A=A\circ a=A,\ \forall a\in A.\eqno(3)$$  
A semihypergroup satisfying $(3)$ is called {\it hypergroup}. 
The factor of an $H_v$-semigroup (or of a semihypergroup) modulo a strongly regular 
equivalence relation is a semigroup and $H_v$-group (or of a hypergroup) modulo a strongly regular 
equivalence relation is a group.
As a matter of fact, the equalities $(3)$ define on $(A,\circ)$ the binary multioperations
$/,\backslash:A\times A\rightarrow P^*(A),$ 
$$
b/a=\{x\in A\mid b\in x\circ a\},\ a\backslash b=\{x\in A\mid b\in
a\circ x\},\eqno(4)
$$ the hypergroups and the $H_v$-groups can be defined using identities 
as in Proposition 1 and Corollary 1 \cite{PePu08}, but $/,\backslash$ need not
be involved in the characterization of their equivalence relations for which
the factor hypergroupoid is a group since   
$E_{ua}(A,\circ)=E_{ua}(A,\circ,/,\backslash)$ (see \cite[Remark 15]{PePu08}).
\end{ex}

\section{The poset $(E_{ua}(\mathfrak A),\subseteq)$ and fundamental relations}\label{s3}

Remind that $(E_{ua}(\mathfrak A),\subseteq)$ is an
algebraic closure system on $A\times A$ (see \cite[Lemma
4.2]{PePu06}), hence $(E_{ua}(\mathfrak A),\subseteq)$
is a complete lattice. We denote by $\alpha^{\mathfrak A}$ (or by
$\alpha$ when it is obvious in which multialgebra we are working)
the corresponding closure operator. The smallest relation from
$E_{ua}(\mathfrak A)$ which contains a relation $R\subseteq A\times
A$ is
$$\alpha^{\mathfrak A}(R)=\alpha(R)=\bigcap\{\rho\in E_{ua}(\mathfrak A)\mid
R\subseteq
\rho\}.$$
The smallest element of $E_{ua}(\mathfrak A)$ is called
{\it the fundamental relation of} $\mathfrak A$ and it is denoted by
$\alpha_{\mathfrak A}^*$. Characterizations of this
relation are presented in \cite[Theorem 2 and Corollary
11]{PePu08}, and \cite[Corollary 4.5]{PePu06}.

\medskip

Let $\mathfrak A=(A,(f_\gamma)_{\gamma<o(\tau)})$ be a multialgebra,
$E(A)$ the set of the equivalence relations of $A$, $\rho\in E(A)$, and $\pi_{\rho}$, 
$\pi_{\alpha^*_{\mathfrak A}}$ the canonical projections of the factor multialgebras
$\mathfrak A/\rho$ and $\mathfrak A/\alpha^*_{\mathfrak A}$, respectively. If 
$\alpha_{\mathfrak A}^*\subseteq\rho$, according to \cite[Theorem 5.5(ii)]{EKM}, there exists a unique multialgebra homomorphism
$\varphi$ which makes the following diagram commutative
$$\xymatrix{ A \ar[r]^{\pi_{\alpha^*_{\mathfrak A}}\ \ \ }\ar[dr]_{\pi_{\rho} } &  A/\alpha^*_{\mathfrak A}\ar@{.>}[d]^{\varphi}  \\
\ &  A/\rho  \ .}$$ 
If we denote by $\rho/\alpha^*_{\mathfrak A}$ the kernel of $\varphi$, then $A/\rho$ is a universal
algebra if and only if $\varphi$ satisfies $(1')$ for any $\gamma<o(\tau)$. According to 
\cite[Theorem 1]{Pi67} and Example \ref{congA}, this happens exactly when 
$\rho/\alpha^*_{\mathfrak A}$ is a congruence relation of the universal algebra
$\mathfrak A/\alpha^*_{\mathfrak A}$. 
Thus, the relations from $E_{ua}(\mathfrak A)$
are also characterized by the following lemma.

\begin{lem}\label{euacong}
If $\mathfrak A=(A,(f_\gamma)_{\gamma<o(\tau)})$ be a multialgebra, then
$$E_{ua}(\mathfrak A)=\{\rho\in E(A)\mid \alpha_{\mathfrak A}^*\subseteq\rho,\ 
\rho/\alpha^*_{\mathfrak A}\in C(\mathfrak A/\alpha^*_{\mathfrak A})\}.$$
\end{lem}

\begin{thm}{\rm\cite[Theorem 4.4]{PePu06}}\label{qr}
Let $\mathbf q, \mathbf r$  be $n$-ary terms $(n\in\mathbb N^*)$, and let
$\mathfrak A$ be a multialgebra of type $\tau$. Let $\alpha_{\mathbf{qr}}\subseteq 
A\times A$ be the relation defined
as follows: $x\alpha_{\mathbf{qr}} y$ if and only if there exist a unary
polynomial function $p\in {\rm Pol}^A_1(\mathfrak P^*(\mathfrak A))$ and
$a_0,\ldots,a_{n-1}\in A$ such that
\begin{align*}
&x\in p(q(a_0,\ldots,a_{n-1})),\ y\in p(r(a_0,\ldots,a_{n-1}))\ \hbox{or}\\
&y\in p(q(a_0,\ldots,a_{n-1})),\ x\in p(r(a_0,\ldots,a_{n-1})).
\end{align*}
The transitive closure $\alpha^*_{\mathbf{qr}}$ of $\alpha_{\mathbf{qr}}$ is the
smallest equivalence relation of $\mathfrak A$ such that the factor
multialgebra is a universal algebra satisfying the identity
$\mathbf q=\mathbf r$.
\end{thm}

Let $I$ be a set and, for any $i\in I$, let $\mathbf q_i,\mathbf r_i$ be $m_i$-ary
terms of type $\tau$ ($m_i\in\mathbb N^*$) and consider the following set of
identities
$$\mathcal I=\{\mathbf q_i=\mathbf r_i\mid i\in I\}.$$ For a multialgebra $\mathfrak A$, we denote
$$R_{\mathcal I}=\bigcup\{q_i(a_0,\ldots,a_{m_i-1})\times
r_i(a_0,\ldots,a_{m_i-1})\mid a_0,\ldots
,a_{m_i-1}\in A,\ i\in I\}.$$ In particular, if $|I|=1$  
and the considered identity is $\mathbf q=\mathbf r$ we will write
$R_{\mathbf{qr}}$ instead of $R_{\mathcal I}$.

In \cite[Remark 13]{PePu06} it was noticed that the smallest equivalence relation of a
multialgebra $\mathfrak A$ of type $\tau$ for which the factor multialgebra is a universal algebra
satisfying an identity $\mathbf q=\mathbf r$ is $\alpha(R_{\mathbf{qr}}).$ We also have:

\begin{lem}\label{l0}
The smallest equivalence relation of a
multialgebra $\mathfrak A$ of type $\tau$
for which the factor multialgebra is a universal algebra
satisfying all the identities from $\mathcal I$ is
the relation $\alpha(R_{\mathcal I})$. In particular, $\alpha^*_{\mathbf{qr}}=\alpha(R_{\mathbf{qr}}).$
\end{lem}

\begin{pf}
Let $\rho\in E_{ua}(\mathfrak A)$. If an identity $\mathbf q_i=\mathbf r_i$
from $\mathcal I$ is satisfied in $\mathfrak A/\rho$ and 
$a_0,\ldots ,a_{m_i-1}\in A$, we deduce from $(2)$ that for any
$(x,y)\in q_i(a_0,\ldots,a_{m_i-1})\times r_i(a_0,\ldots,a_{m_i-1}),$
$$\rho\langle x\rangle=q_i(\rho\langle a_0\rangle,\ldots,\rho\langle a_{m_i-1}\rangle)=
r_i(\rho\langle a_0\rangle,\ldots,\rho\langle a_{m_i-1}\rangle)=\rho\langle y\rangle,$$
hence $(x,y)\in\rho$. Thus, if all the identities from $\mathcal I$
are valid in $\mathfrak A/\rho$, then $R_{\mathcal I}\subseteq\rho$. 
Conversely, if $R_{\mathcal I}\subseteq\rho$,
using again $(2)$, we easily deduce that all the identities from $\mathcal I$
hold in $\mathfrak A/\rho$. 
It follows that any equivalence relation from $E_{ua}(\mathfrak A)$ for which the factor
multialgebra is a universal algebra satisfying the identities from $\mathcal I$ must contain 
$R_{\mathcal I}$. Thus, the smallest such relation is $\alpha(R_{\mathcal I})$.
\end{pf}

In the proof of the above lemma we showed that if all the identities from $\mathcal I$ are satisfied on 
$\mathfrak A/\rho$ then $R_{\mathcal I}\subseteq\rho$. Applying the closure operator $\alpha$ to this 
inclusion, we obtain:

\begin{cor}\label{c1}
If $\rho\in E_{ua}(\mathfrak A)$ and all the identities from $\mathcal I$ are satisfied on 
$\mathfrak A/\rho$ then $\alpha(R_{\mathcal I})\subseteq\rho\,.$
\end{cor}

\begin{lem}
In the complete lattice $(E_{ua}(\mathfrak
A),\subseteq)$ we have $$\alpha(R_{\mathcal I})=\bigvee_{i\in
I}\alpha^*_{\mathbf q_i\mathbf r_i}.$$
\end{lem}

\begin{pf} If $I=\emptyset$, the equality holds since 
$R^{ }_{\mathcal I}=\emptyset$ and $\alpha(\emptyset)$ is the 
fundamental relation of $\mathfrak A$, which is the 
infimum of $E_{ua}(\mathfrak A),$ hence the supremum of any empty family 
of relations from $E_{ua}(\mathfrak A).$ 
If $I\neq\emptyset$ then, according to Lemma \ref{l0} and to the properties 
of the closure operator, we have  
$$\alpha(R_{\mathcal I})=\alpha\left(\bigcup_{i\in I} R_{\mathbf q_i\mathbf r_i}\right)
=\bigvee_{i\in I}\alpha(R_{\mathbf q_i\mathbf r_i})=\bigvee_{i\in
I}\alpha^*_{\mathbf q_i\mathbf r_i}.$$
\end{pf}

\begin{rmk}\label{r3}
From Remark \ref{r2} it follows that for
any set of identites $\mathcal I$ satisfied at least in a weak manner on the multialgebra 
$\mathfrak A$, $\alpha(R_{\mathcal I})$ is the smallest element from $E_{ua}(\mathfrak
A)$, i.e. $\alpha_{\mathfrak A}^*=\alpha(R_{\mathcal I}).$ In particular, 
for any variable $\mathbf x$, $$\alpha_{\mathfrak
A}^*=\alpha^*_{\mathbf{xx}}=\alpha(R_{\mathbf{xx}}),$$ hence, 
if $I=\emptyset$, then $\alpha(R_{\mathcal I})=\alpha_{\mathfrak A}^*=\alpha(R_{\mathbf{xx}})$.
Thus, we can consider $I\neq\emptyset$ without loosing the fundamental
relation from our study.   
\end{rmk}

\begin{thm}\label{qrI}
Let $I\neq\emptyset$ be a set and for any $i\in I$, let $\mathbf q_i,\mathbf r_i$ be $m_i$-ary
terms of type $\tau$, let $\mathcal I=\{\mathbf q_i=\mathbf r_i\mid i\in I\}$, and let 
$\mathfrak A=(A,(f_\gamma)_{\gamma<o(\tau)})$ be a multialgebra of type $\tau$. Let
$\alpha_{\mathcal I}\subseteq A\times A$ be the relation defined as
follows:
\begin{align*}
x&\alpha_{\mathcal I} y\Leftrightarrow \exists i\in I,\ \exists
p_i\in {\rm Pol}^A_1(\mathfrak P^*(\mathfrak A)),\ \exists
a^i_0,\ldots,a^i_{m_i-1}\in
A:\\
&x\in p_i(q_i(a^i_0,\ldots,a^i_{m_i-1})),\ y\in p_i(r_i(a^i_0,\ldots,a^i_{m_i-1}))\
\hbox{or}\\
&y\in p_i(q_i(a^i_0,\ldots,a^i_{m_i-1})),\ x\in p_i(r_i(a^i_0,\ldots,a^i_{m_i-1})).
\end{align*}
The transitive closure $\alpha^*_{\mathcal I}$ of
$\alpha_{\mathcal I}$ is the smallest equivalence relation on $\mathfrak A$
for which the factor multialgebra is a universal algebra
satisfying all the identities from $\mathcal I$, i.e.
$$\alpha^*_{\mathcal I}=\alpha(R_{\mathcal I}).$$
\end{thm}

\begin{pf}
From Lemma \ref{euacong}, it follows that for any relations $\rho_i\in E_{ua}(\mathfrak A)$, $i\in I$,
their supremum $\sup_{i\in I}\rho_i$ in $(E(A),\subseteq)$ is 
$\bigvee_{i\in I}\rho_i$ (i.e. their supremum in $(E_{ua}(\mathfrak A),\subseteq)$).
Indeed, $\alpha^*_{\mathfrak A}\subseteq\sup_{i\in I}\rho_i$,
$(\sup_{i\in I}\rho_i)/\alpha^*_{\mathfrak A}=\sup_{i\in I}\rho_i/\alpha^*_{\mathfrak A}$ and,
since the congruence lattice of a universal algebra is a complete sublattice of the equivalence lattice
of its universe, $\sup_{i\in I}\rho_i/\alpha^*_{\mathfrak A}\in C(\mathfrak A/\alpha^*_{\mathfrak A}).$

Thus $$\sup_{i\in I}\alpha^*_{\mathbf q_i\mathbf r_i}=\bigvee_{i\in I}\alpha^*_{\mathbf q_i\mathbf r_i}
=\alpha(R_{\mathcal I}).$$
But $\alpha_{\mathbf q_i\mathbf r_i}$ ($i\in I$) are symmetric and reflexive relations, hence 
$\sup_{i\in I}\alpha^*_{\mathbf q_i\mathbf r_i}$ is the transitive closure of the union 
$\bigcup_{i\in I}\alpha_{\mathbf q_i\mathbf r_i}=\alpha_{\mathcal I}$, which is exactly
$\alpha^*_{\mathcal I}$.
\end{pf}

Using Corollary \ref{c1} and the notations from the above theorem, we have:

\begin{cor}\label{c3}
If $\mathcal J$ is a set of identities satisfied 
in the universal algebra $\mathfrak A/\alpha^*_{\mathcal I}$ then 
$\alpha^*_{\mathcal J}\subseteq\alpha^*_{\mathcal I}.$
\end{cor}

Since $\alpha^*_{\mathcal I}$ is the smallest equivalence relation on $\mathfrak A$
for which the factor multialgebra is a universal algebra
satisfying all the identities from $\mathcal I$ we also have:

\begin{cor}\label{c4}
If $\mathcal J$ is a set of identities satisfied 
in the universal algebra $\mathfrak A/\alpha^*_{\mathcal I}$ then 
$\alpha^*_{\mathcal I\cup\mathcal J}=\alpha^*_{\mathcal I}.$
\end{cor}

In the next sections we will refer to the relation $\alpha^*_{\mathcal I}$ of $\mathfrak A$ as
the $\mathcal I$-{\it fundamental relation of} $\mathfrak A$ and to 
$\mathfrak A/\alpha^*_{\mathcal I}$ as the $\mathcal I$-{\it fundamental algebra of} $\mathfrak A$.

\section{A new characterization of the $\alpha^*$-relation of a hyperring}\label{sectionhr}

In the literature there are many examples of relations as those from Theorem \ref{qrI} defined
on particular hyperstructures. In many cases, the authors want to obtain at most one identity in the factor (multi)algebra.
Such examples are the fundamental relations which appear in \cite{LeDa08} and the relation $\gamma^*$ from \cite{Fr02}.
These situations are rather related to \cite{PePu06}.  
The $\alpha^*$-relations from \cite{DaVo} are very illustrative for the case when the set $\mathcal I$ 
contains more than one identity. This is why we will use \cite{DaVo} 
to exemplify our construction and also to give some hints on how our general results can be helpful in 
approaching particular multialgebra cases. 

\smallskip 

A set $A$ with two binary multioperations $+$ and $\cdot$ is
called {\it $H_v$-ring} if $(A,+)$ is an $H_v$-group, $(A,\cdot)$
is an $H_v$-semigroup, and for any $a,b,c\in A$
$$a(b+c)\cap(ab+ac)\neq\emptyset\ \hbox{and}\
(b+c)a\cap(ba+ca)\neq\emptyset.$$
An $H_v$-ring $(A,+,\cdot)$ is called {\it hyperring} if $(A,+)$ is hypergroup, $(A,\cdot)$
is semihypergroup, and for any $a,b,c\in A$
$$a(b+c)=(ab+ac)\ \hbox{and}\ (b+c)a=(ba+ca).$$

\begin{lem}{\rm\cite[Corollary 9]{PePu08}}
Let $(A,+,\cdot)$ be an $H_v$-ring. An equivalence relation $\rho$
on $A$ is in $E_{ua}(A,+,\cdot)$ if and only if it is strongly
regular both on $(A,+)$ and $(A,\cdot)$. In this case, the factor
multialgebra $(A/\rho,+,\cdot)$ is a nearring. More exactly, $+$
needs to be commutative for $(A/\rho,+,\cdot)$ to be a ring.
\end{lem}

\begin{rmk}
According to Example \ref{exhypslash}, the hyperrings ($H_v$-rings) can be seen as multialgebras
$(A,+,/,\backslash,\cdot)$ with two binary (weak) associative multioperations $+$, $\cdot$, and two
binary multioperations $/,\backslash$ satisfying $(4)$. Since {\small
$$E_{ua}(A,+,/,\backslash,\cdot)=E_{ua}(A,+,/,\backslash)\cap E_{ua}(A,\cdot)=
E_{ua}(A,+)\cap E_{ua}(A,\cdot)=E_{ua}(A,+,\cdot),$$
 }we do not have to use the multioperations $/$ and $\backslash$ in the construction of the 
polynomial functions which appear in the characterization or in the construction of the 
relations from $E_{ua}(A,+,/,\backslash,\cdot)$.
\end{rmk}

Applying Theorem \ref{qrI} to an $H_v$-ring $(A,+,\cdot)$ and 
$$\mathcal I=\{\mathbf x_0+\mathbf x_1=\mathbf x_1+\mathbf x_0,\ 
\mathbf x_0\cdot \mathbf x_1=\mathbf x_1\cdot \mathbf x_0\},$$  one obtains:

\begin{prop}\label{p2}
The smallest equivalence relation of an $H_v$-ring $(A,+,\cdot)$ for which the factor multialgebra  
is a commutative ring is the transitive closure $\alpha^*_{\mathcal I}$ of the relation $\alpha_{\mathcal I}$ 
defined by
\begin{align*}
x\,\alpha_{\mathcal I}\, y\ \Leftrightarrow\  &\exists
p_0,p_1\in {\rm Pol}^A_1(\mathfrak P^*(A,+,\cdot)),\ \exists
a^0_0,a^0_1,a^1_0,a^1_1\in
A:\\
&x\in p_0(a^0_0+a^0_1),\ y\in p_0(a^0_1+a^0_0)\
\hbox{or}\\
&x\in p_1(a^1_0\cdot a^1_1),\ y\in p_1(a^1_1\cdot a^1_0).
\end{align*}
\end{prop}

From the equivalence of the definitions of a
polynomial function from \cite{BuSa} and \cite{Gr79} it follows:

\begin{lem}\label{l3}
For any unary polynomial function $p\in {\rm Pol}^A_1(\mathfrak
P^*(A,+,\cdot))$, there exist $m\in\mathbb N$, $m\geq 1$, $b_1,\ldots, b_{m-1}\in A$ and a term
function $p'\in {\rm Clo}_m(\mathfrak P^*(A,+,\cdot))$ such that for any $A_0\in P^*(A)$,
$$p(A_0)=p'(A_0, b_1,\ldots ,b_{m-1}).$$ 
\end{lem}

Using Lemma \ref{l3}, we can replace the unary polynomial functions from Proposition \ref{p2}
by term functions. In hyperrings, the distributivity of $\cdot$ with respect to $+$ allows us to use in 
the characterization of $\alpha^*_{\mathcal I}$ only the term functions $p'$ which are induced by 
the terms $\mathbf p'$ which are sums of products of a finite number of variables and to rewrite 
the relation $\alpha^*_{\mathcal I}$ as in \cite{DaVo}. Thus, for the next part of this section, we will 
consider $(A,+,\cdot)$ to be a hyperring. Since the closure operator $\alpha$ will not explicitly appear in 
this section, we think it is alright to use the notations from \cite{DaVo}.

\begin{dfn}\cite[Definition 1]{DaVo}\label{d1}
For a hyperring $(A,+,\cdot)$, one defines the relation $\alpha^*$ as the transitive closure of
the relation
\begin{align*}
x\alpha y \Leftrightarrow &\exists n,k_1,\ldots,k_n\in\mathbb N^*,\, \exists 
\sigma\in S_n,\, \exists x_{i1},\ldots,x_{ik_i}\in A,\, \exists \sigma_i\in S_{k_i} 
(i=1,\ldots,n):\\
&x\in\sum_{i=1}^n\left(\prod_{j=1}^{k_i}x_{ij}\right) \ \hbox{and}\
y\in \sum_{i=1}^n A_{\sigma(i)},\ \hbox{where}\ A_i=\prod_{j=1}^{k_i}x_{i\sigma_i(j)},
\end{align*}
i.e. $y$ is in a sum of products obtained from $\sum_{i=1}^n\left(\prod_{j=1}^{k_i}x_{ij}\right)$
by permuting (according to $\sigma$) the terms of the sum and the factors in each product $A_i$ (as
the permutations $\sigma_i$ indicate).  
\end{dfn}

The following result shows that the relation $\alpha^*$ defined above is in $E_{ua}(A,+,\cdot)$. 

\begin{lem}{\rm\cite[Lemma 2]{DaVo}}\label{l4}
 The relation $\alpha^*$ is a strongly regular equivalence relation both on $(A,+)$ and $(A,\cdot)$.
\end{lem}

From the above lemma we have $\alpha^*_{\mathcal I}\subseteq\alpha^*$ and from
Proposition \ref{p2} we deduce that $(A/\alpha^*_{\mathcal I},+,\cdot)$ is a commutative ring.
It results that if we take any sum of products of elements from $A/\alpha^*_{\mathcal I}$ and permute the 
factors of the products and the terms of the sums (including the ``degenerated'' cases when the sum has only
one term, or some products have only one factor, or both), we obtain the same element. 
If $\mathcal J$ is the set of identities which determine all these identifications then
$\alpha\subseteq R_{\mathcal J}\subseteq\alpha^*_{\mathcal J}$, hence $\alpha^*\subseteq\alpha^*_{\mathcal J}.$ 
Applying Corollary \ref{c3}, we have $\alpha^*_{\mathcal J}\subseteq\alpha^*_{\mathcal I}$, so
$\alpha^*\subseteq\alpha^*_{\mathcal I}$, and, consequently, $$\alpha^*=\alpha^*_{\mathcal I}.$$
Thus, using Lemma \ref{l4}, Proposition \ref{p2}, and Corollary \ref{c3}, one can avoid \cite[Theorem 3]{DaVo}
(and the computations from the its proof) in order to obtain:

\begin{thm}{\rm\cite[Theorem 4]{DaVo}}\label{t4DaVo}
 The relation $\alpha^*$ is the smallest equivalence relation on $(A,+,\cdot)$ such that the 
factor multialgebra is a commutative ring.
\end{thm}

From Remark \ref{r3} we immediately deduce:

\begin{thm}{\rm\cite[Theorem 13]{DaVo}}
If the multioperations of the hyperring $(A,+,\cdot)$ are weak commutative, then $\alpha^*$
is equal to the fundamental relation of $(A,+,\cdot)$.
\end{thm}

The following result is an improvement for the form of $\alpha^*_{\mathcal I}$: 

\begin{thm}\label{t3}
The smallest equivalence relation of a hyperring $(A,+,\cdot)$ for which the factor multialgebra  
is a commutative ring is the transitive closure $\alpha^*_0$ of the relation $\alpha_0$ defined 
as the union of the Cartesian products of (finite) sums of (finite) products of elements from $A$ 
(including the cases when the sums have only one term or/and the products have only one factor),
$$t(a_0, a_1,\ldots ,a_{m-1})\times t'(a_0, a_1,\ldots ,a_{m-1}),$$
such that $t'(a_0, a_1,\ldots ,a_{m-1})$ is obtained from $t(a_0, a_1,\ldots ,a_{m-1})$ either by permuting
two consecutive factors in a product or by permuting two consecutive products from the sum. 
\end{thm}

\begin{pf}
We can prove this result noticing that 
$\alpha^*_0$ and the transitive closure $\alpha^*$ of the relation $\alpha$ from Definition \ref{d1} 
are the same. Knowing how the transitive closure of a relation works and the fact that all the transpositions 
of the form $(i-1,\,i)$, $i\in\{2,\ldots,n\}$, generate the group $(S_n,\cdot)$ of the permutations of the set 
$\{1,\ldots,n\}$ it is an easy exercise to show that the relation $\alpha$ 
is included in the transitive closure of $\alpha_0$. Also, it is  obvious that $\alpha_0\subseteq\alpha$. Thus 
$\alpha^*=\alpha^*_0.$
\end{pf}
 
If $\cdot$ is an operation, the products from any $t(a_0, a_1,\ldots ,a_{m-1})$ are singletons. If the operation $\cdot$ is
commutative, the permutation of two consecutive factors in such a product does not change the product, so
these (one-element set) products are the same in $t(a_0, a_1,\ldots ,a_{m-1})$ and $t'(a_0, a_1,\ldots ,a_{m-1})$.
Hence, any $t$ from the previous theorem is a sum of elements from $A$ and any $t'$ is either $t$, or it is obtained from $t$ 
by permuting two consecutive terms of this sum. Thus, in this case, $\alpha_0^*$ is  the transitive closure of the relation
\begin{align*}
x\,\alpha'_0 \,y\ \Leftrightarrow\ &\exists n\in\mathbb N^*,\ \exists z_1,\ldots,z_n\in A,\
\exists i\in\{1,\ldots,n-1\}:\\ &x\in z_1+\cdots+
z_{i-1}+(z_i+ z_{i+1})+ z_{i+2}+\cdots+ z_n,\\ &y\in z_1+\cdots+
z_{i-1}+(z_{i+1}+ z_i)+ z_{i+2}+\cdots+ z_n.
\end{align*}
Using \cite[Example 5]{PePu06}, we have:

\begin{cor}{\rm\cite[Theorem 5]{DaVo}}\label{c5}
For those hyperrings $(A,+,\cdot)$ for which $(A,\cdot)$ is a commutative semigroup, the relation $\alpha^*$ 
is equal to the smallest equivalence relation of $(A,+)$ for which the factor multialgebra is a commutative group.
\end{cor}

\section{On the directed colimits of fundamental algebras}\label{dirlim}

The multialgebras of type $\tau$, multialgebra homomorphisms and 
mapping composition determine a category. We denote it by $\mathbf{Malg}(\tau)$.
In some previous papers, Pelea showed that the factorization modulo the fundamental relation
of a multialgebra determines a functor from $\mathbf{Malg}(\tau)$ to  the category
$\mathbf{Alg}(\tau)$ of the universal algebra of type $\tau$ (see, for instance, 
\cite[Theorem 1]{Pe05}). He also investigated if and when this functor preserves
different category theoretical constructions of multialgebras (see \cite{Pe04, Pe05, Pe06}).  
Recently, the problem was revisited in the case of hyperrings, commutative fundamental relations and finite (direct) 
products (see, \cite{DaVo, MiDa}). The general approach from Section \ref{s3}
provides an easy way to solve the problem in the general case of multialgebras, 
$\mathcal I$-fundamental relations and directed colimits. 

\smallskip

Any variety $\mathcal V$ of universal algebras of type $\tau$ can be seen as a (full) subcategory of the category
$\mathbf{Alg}(\tau)$. 

\begin{thm}\label{refl}
Any variety $\mathcal V$ of universal algebras of type $\tau$ is a reflective subcategory of $\mathbf{Malg}(\tau)$.
\end{thm}

\begin{pf} 
Let $\mathcal I$ be a set of identities such that an algebra is in $\mathcal V$ if and only if it satisfies all
the identities from $\mathcal I$, let $\mathfrak A$ be a multialgebra and let $\pi_{\alpha^*_{\mathcal I}}$ 
be  the canonical projection of $\mathfrak A$ determined by the relation $\alpha^*_{\mathcal I}$. 
Let $\mathfrak B$ be a universal algebra from $\mathcal V$ and let $h:A\to B$ be a multialgebra homomorphism. Since
$\mathfrak B$ is a universal algebra, $h$ is an ideal (or strong) homorphism, hence $h(A)$ is a subalgebra of $\mathfrak B$ so,
$h(\mathfrak A)\in \mathcal V$. From \cite[Theorem 1]{Pi67}, it follows that $h(\mathfrak A)\simeq \mathfrak A/\ker h,$
thus $\ker h\in E_{ua}(\mathfrak A)$ and $\alpha^*_{\mathcal I}\subseteq \ker h$. According to \cite[Theorem 5.5]{EKM},
there exists a unique universal algebra homomorphism $\overline h$ which makes the following diagram
$$\xymatrix{ \mathfrak A \ar[r]^{\pi_{\alpha^*_{\mathcal I}}\ \ \ }\ar[dr]_{h} &  \mathfrak A/\alpha^*_{\mathcal I}\ar@{.>}[d]^{\overline{h}}  \\
\ &  \mathfrak B  }\eqno(5)$$ commutative. This means that $\pi_{\alpha^*_{\mathcal I}}$ is a $\mathcal V$-reflection 
for $\mathfrak A$.
\end{pf}

From \cite[Proposition 4.22]{AHS} one can can easily deduce the following:

\begin{cor}
For a fixed set of identities $\mathcal I$, the factorization modulo 
$\alpha^*_{\mathcal I}$ determines a (covariant) functor $F_{\mathcal I}:\mathbf{Malg}(\tau)\to \mathcal V$. 
\end{cor}

The functor $F_{\mathcal I}$ is defined as follows: for any multialgebra $\mathfrak A$ of type $\tau$, $F_{\mathcal I}
(\mathfrak A)=\mathfrak A/\alpha^*_{\mathcal I}$, and for any multialgebra homomorphism $\mathfrak A\overset{h}{\to}\mathfrak B$ 
from $\mathbf{Malg}(\tau)$, $F_{\mathcal I}(h)$ is the unique algebra homomorphism which makes the following diagram
$$\diagram
            \mathfrak A\rto^{h}\dto_{\pi^{\mathfrak A}_{\alpha^*_{\mathcal I}}} &\mathfrak B\dto^{\pi^{\mathfrak B}_{\alpha^*_{\mathcal I}}}\\
            \mathfrak A/\alpha^*_{\mathcal I}\rto^{F_{\mathcal I}(h)} &\mathfrak B/\alpha^*_{\mathcal I}
\enddiagram$$ commutative (we denoted by $\pi^{\mathfrak A}_{\alpha^*_{\mathcal I}}$ and $\pi^{\mathfrak B}_{\alpha^*_{\mathcal I}}$ 
the canonical projections of $\mathfrak A$ and $\mathfrak B$, respectively, determined by their $\mathcal I$-fundamental relations).

\smallskip

If we take $\mathcal I=\{\mathbf x=\mathbf x\}$ for some variable $\mathbf x$, then $\alpha^*_{\mathcal I}=\alpha^*_{\mathfrak A}$, 
and this functor is the functor $F:\mathbf{Malg}(\tau)\to \mathbf{Alg}(\tau)$ studied in \cite{Pe04, Pe05} and \cite{Pe06}. 
In \cite{Pe06}, Pelea constructed the colimit of a directed diagram $D:(I,\leq)\to \mathbf{Malg}(\tau)$ 
(called direct system of multialgebras in \cite{Pe06}). Theorem \ref{refl} strengthens 
\cite[Theorem 25]{Pe06} since the reflector $F_{\mathcal I}$ is a coadjoint for the inclusion functor 
$\mathcal V\hookrightarrow \mathbf{Malg}(\tau)$. As a matter of fact, $\mathcal V$ is the subcategory of the 
reflexive objects with respect to the counit of the above adjunction (see, for instance, \cite[Section 1]{Mo10}).
But coadjoint functors preserve the existing colimits, thus we have:

\begin{prop}\label{Icolim}
Let $D:(I,\leq)\to \mathbf{Malg}(\tau)$ be a directed diagram. The $\mathcal I$-fundamental algebra of the directed 
colimit $\text{\rm colim} D$ is isomorphic to the (directed) colimit of the composition $(I,\leq)\overset{D}{\to} \mathbf{Malg}(\tau)
\overset{F_{\mathcal I}}{\to}\mathcal V$.
\end{prop}

The $\mathcal I$-fundamental algebra of a multialgebra preserves multialgebra's identities (see
Remark \ref{r2}). Thus, the diagram $(5)$ is commutative if we take 
$$\mathcal I=\{\mathbf x_0+\mathbf x_1=\mathbf x_1+\mathbf x_0,\ 
\mathbf x_0\cdot \mathbf x_1=\mathbf x_1\cdot \mathbf x_0\},$$ $\mathfrak A$ an $H_v$-ring
(or a hyperring) and $\mathfrak B$ a commutative ring. One deduces the following:

\begin{prop}
The category $\mathbf{CRng}$ of commutative rings is a (full) reflective subcategory of the category
$\mathbf H_v\mathbf{Rng}$ of $H_v$-rings (and $H_v$-ring homomorphisms)
and also of the category $\mathbf{HR}$ of hyperrings (and hyperring homomorphisms). 
\end{prop}
 
Also, if $D:(I,\leq)\to \mathbf{Malg}(\tau)$ is a directed diagram and all the multialgebras $D(i)$, $i\in I$, 
satisfy a given (weak or strong) identity, the colimit of $D$ also satisfies the given identity 
(see \cite[Lemmas 29 and 31]{Pe06}). Therefore, the result from Proposition \ref{Icolim} holds for
$H_v$-rings and hyperrings.

\begin{prop}
The commutative-fundamental ring of a directed colimit of $H_v$-rings 
(or hyperrings) is isomorphic to the directed colimit of the corresponding commutative-fundamental rings.
\end{prop}

 \end{document}